\documentclass{article}

\usepackage{proof}
\usepackage{latexsym}
\usepackage{amssymb}

\newtheorem{theorem}{Theorem}
\newtheorem{definition}[theorem]{Definition}
\newtheorem{proposition}[theorem]{Proposition}
\newtheorem{lemma}[theorem]{Lemma}

\title{Provably $\Delta^{0}_{2}$ and weakly descending chains
\thanks{Dedicated to the occasion of Chong Chi Tat's 60th birthday}
}
\author{Toshiyasu Arai
\\
Graduate School of Science,
Chiba University
\\
1-33, Yayoi-cho, Inage-ku,
Chiba, 263-8522, JAPAN
\\
tosarai@faculty.chiba-u.jp
}
\date{}
\begin{document}
\maketitle

\begin{abstract}
In this note we show that a set is provably $\Delta^{0}_{2}$ in the fragment $I\Sigma_{n}$ 
of arithmetic iff
it is $I\Sigma_{n}$-provably in the class $D_{\alpha}$ of $\alpha$-r.e. sets in the Ershov
hierarchy for an $\alpha<_{\varepsilon_{0}}\omega_{1+n}$, where $<_{\varepsilon_{0}}$ denotes a standard $\varepsilon_{0}$-ordering.

In the Appendix it is shown that a limit existence rule $(LimR)$ due to Beklemishev and Visser
becomes stronger when the number of nested applications of the inference rule grows.
\end{abstract}

\section{Introduction}

Thoroughout this paper, we identify a predicate $A$ with its characteristic function
\[
A(x_{1},\ldots,x_{n})=\left\{
\begin{array}{ll}
0 & \mbox{{\rm if }} A(x_{1},\ldots,x_{n})
\\
1 & \mbox{{\rm otherwise}}
\end{array}
\right.
\]

Natural numbers $c$ are identified with the sets $\{n\in\mathbb{N}: n<c\}$.

The following Limit Lemma due to Shoenfield is a classic in computability theory.

\begin{theorem}\label{th:limitlemma}{\rm (Limit Lemma)}\\
A set $A$ of natural numbers is $\Delta^{0}_{2}$ iff there is a binary (primitive) recursive predicate $f:\omega\times\omega\to 2=\{0,1\}$ such that
\[
\forall c[\lim_{w\to\infty}f(c,w)=A(c)]
.\]
\end{theorem}

Moreover the theorem is provable uniformly in $B\Sigma^{0}_{1}\subseteq I\Sigma^{0}_{1}$, 
cf. \cite{H-P}, pp. 89-91.
Let us call the predicate $f$ a {\it witnessing predicate\/} for $A\in\Delta^{0}_{2}$.

In this paper we address a problem asking what can we say about the rate of convergences
of the predicate $f$ under the assumption that the set $A$ is {\it provably\/} $\Delta^{0}_{2}$ 
in a formal (sound) theory T?

This is a problem on a hierarchy.
The class of $\Delta^{0}_{2}$-sets is classified in the Ershov hierarchy, \cite{Ershov2}.
A recent article \cite{alcSYY} due to F. Stephan, Y. Yang and L. Yu is a readable contribution to the hierarchy, to which we refer as a standard text.

The $\alpha$-th level of the Ershov hierarchy is denoted $D_{\alpha}$ for notations $\alpha$ 
of constructive ordinals, and a set in $D_{\alpha}$ is said to be an $\alpha$-r.e. set.

It is known, as usual in hierarchic problems indexed by constructive ordinals, that $D_{\alpha}$ depends
 heavily on notations $\alpha$, i.e., the order type of $\alpha$ does not determine the set $D_{\alpha}$.
By reason of this dependency let us fix a standard elementary recursive well ordering $<_{\alpha}$ of type $\alpha$.
I don't want to discuss here what is a 'standard ordering' or a 'natural well ordering'.
We assume that {\sf EA}=$\mbox{I}\Delta^{0}_{0}+exp$, Elementary Recursive Arithmetic, proves some algebraic
facts on the ordering $<_{\alpha}$.
For the case $\alpha=\varepsilon_{0}$, what we need on $<_{\varepsilon_{0}}$ can be found in, e.g., \cite{Sommer}.

In what follows let us drop the subscript $\alpha$ in $<_{\alpha}$ when no confusion likely occurs.

\begin{definition}\label{df:alcSYY4.1}{\rm (Stephan-Yang-Yu \cite{alcSYY})}\\
{\rm Let} $K\in dom(<)${\rm , the domain of the order} $<$.

{\rm A set} $A$ {\rm of natural numbers is} $K$-r.e.  with respect to $<$ {\rm iff there exist
 a binary}
 recursive {\rm predicate} $f${\rm , and a}
 recursive {\rm function}
$h:\omega \times\omega\to K=\{\beta\in dom(<): \beta<K\}$ such that 
\begin{enumerate}
\item 
\begin{equation}\label{th:delta2pra.1SYY}
\mbox{{\rm (weakly descending)} }  K>h(c,w)\geq h(c,w+1)
\end{equation}
\item  
\begin{equation}\label{th:delta2pra.2SYY}
\mbox{{\rm (lowering)} }  f(c,w)\neq f(c,w+1) \to h(c,w)>h(c,w+1)
\end{equation}
\item 
\begin{equation}\label{th:delta2pra.3SYY}
\forall c[\lim_{w\to\infty}f(c,w)=A(c)]
\end{equation}
\end{enumerate}

{\rm Roughly speaking, a set is} $K${\rm -r.e. if the convergence of its witnessing predicate follows
from the fact that weakly decreasing functions in} $K$ {\rm have to be constant eventually.}
\end{definition}

Now suppose that we have a proof-theoretic analysis of a formal (and sound) theory T, e.g., 
a cut-elimination through a transfinite induction along a standard well ordering $<$.
It, then, turns out that
$A$ is provably $\Delta^{0}_{2}$ in T iff T proves the fact that $A\in D_{K}$ with respect to $<$ for a $K\in dom(<)$.

Though, in this paper, we restrict our attention to $\mbox{T}=I\Sigma^{0}_{n}$ of fragments of first order arithmetic
as a concrete example, where the order $<$ denotes a standard well ordering of type $\varepsilon_{0}$,
it is easy to see that our proof works also for stronger theories, e.g., second order arithmetic
$\Pi^{1}_{1}\mbox{-CA}_{0}$ and fragments of set theories.

In Section \ref{sect:ISigman} it is shown that for each $n\geq 1$,
 a set is provably $\Delta^{0}_{2}$ in the fragment $I\Sigma_{n}$ iff
it is $I\Sigma_{n}$-provably in the class $D_{\alpha}$ for an $\alpha<_{\varepsilon_{0}}\omega_{1+n}$ (Theorem \ref{th:delta2Isigman}).

Also 
any provably $\Sigma^{0}_{2}$-function has a Skolem function $F(c)= \lim_{w\to\infty}f(c,w)$ as limits of an $f$,
whose convergence is ensured by weakly descending chains of ordinals
(Theorem \ref{th:sigma2witness}).
Moreover the 2-consistency $\mbox{RFN}_{\Pi^{0}_{3}}(I\Sigma^{0}_{n})$
is seen to be equivalent over Primitive Recursive Arithmetic {\sf PRA} to the fact that every primitive recursive 
weakly descending chain of ordinals$<\omega_{1+n}$ has a limit(Theorem \ref{th:2con}).

In Section \ref{sect:PRA} it is shown that a set is provably $\Delta^{0}_{2}$ in
Elementary Recursive Arithmetic  {\sf EA}
 iff
it is {\sf EA}-provably in the class $D_{n}$ of a finite level (Theorem \ref{th:delta2pra}).
Our proof seems to be a neat application of the Herbrand's theorem.

The Appendix \ref{appendix} contains another application of Herbrand's theorem.
We consider, over {\sf EA},
 an inference rule $(LimR)$ in \cite{Lev},
which concludes the convergence of
 an elementary recursive series $\{h(n)\}_{n}$
under the assumption that the series is weakly decreasing almost all $n$.
Note that $(LimR)$ is an inference rule, and not an axiom(sentence).

On the other side, let $L\Sigma_{1}^{-(k)}$ denote the schema in \cite{KPD}, saying that 
any non-empty $\Sigma^{0}_{1}$ $k$-ary predicate has the least tuple, which is least with respect to the lexicographic ordering on $\mathbb{N}^{k}$.

It is shown that $L\Sigma_{1}^{-(k)}$ is equivalent to the $k$-nested applications of $(LimR)$.
In \cite{KPD}, Corollary 2.11 it was shown that $\{L\Sigma_{1}^{-(k)}\}_{k}$ forms a proper hierarchy, i.e.,
$L\Sigma_{1}^{-(k+1)}\vdash\mbox{Con}(L\Sigma_{1}^{-(k)})$.
Hence we conclude that a $(k+1)$-nested application of $(LimR)$ proves the consistency
of the $k$-nested applications of $(LimR)$.

\section{Provably $\Delta^{0}_{2}$ in $I\Sigma^{0}_{n}$}\label{sect:ISigman}

Let {\sf LEA} [{\sf EA}] denote the Lower Elementary Recursive Arithmetic [Elementary Recursive Arithmetic], which is a first-order theory
in the language having function constants for each code(algorithm) of 
lower elementary recursive function [function constants for each code of  elementary recursive function], resp.
Cf. \cite{Rose} and \cite{Skolem} for these classes of subrecursive functions.
Induction schema is restricted to quantifier-free formulas in the language.
The axioms of the theories {\sf LEA}, {\sf EA} are purely universal ones.

Let $I\Sigma^{0}_{n}$ denote the fragment of arithmetic, which is a first-order theory
in the language of {\sf LEA}, and Induction schema is restricted to $\Sigma^{0}_{n}$ formulas.
Here a $\Sigma^{0}_{0}$ formula is a quantifier-free formula.
$I\Sigma^{0}_{0}$ is another name for {\sf LEA}.

Let $<_{\varepsilon_{0}}$ denote a standard $\varepsilon_{0}$-ordering.
We assume that {\sf EA} proves some algebraic
facts on the ordering $<_{\varepsilon_{0}}$.
What we need on $<_{\varepsilon_{0}}$ can be found in, e.g., \cite{Sommer}.

In what follows let us drop the subscript $\varepsilon_{0}$ in $<_{\varepsilon_{0}}$ when no confusion likely occurs.

For a class $\Phi$ of formulas and an ordinal $\alpha$ let $TI(\Phi,\alpha)$ denote the schema of
transfinite induction up to $\alpha$ and applied to a formula $\varphi\in\Phi$:
\[
\forall \beta[ \forall\gamma<\beta \varphi(\gamma) \to \varphi(\beta) ] \to \forall \beta<\alpha\varphi(\beta)
.\]

Let
\[
\omega_{0}:=1, \: \omega_{1+n}:=\omega^{\omega_{n}}
.\]

Here is a folklore result on provability of the restricted transfinite induction schemata in 
fragments of arithmetic.

\begin{theorem}\label{th:Sommer}{\rm (See, e.g., \cite{Sommer})}\\
For each $n\geq 0$,
$I\Sigma^{0}_{n} \vdash TI(\Pi^{0}_{1},\alpha)$ iff $\alpha<\omega_{1+n}$.
\end{theorem}

The following Theorem \ref{th:delta2Isigman} states that for positive integers $n$, 
a set is provably $\Delta^{0}_{2}$ in $I\Sigma^{0}_{n}$
 iff
it is $I\Sigma^{0}_{n}$-provably in the class $D_{\alpha}$ of $\alpha$-r.e. sets in the Ershov hierarchy for an $\alpha<_{\varepsilon_{0}}\omega_{1+n}$.
Moreover (weakly descending) and (lowering) are provable in {\sf EA}.

\begin{theorem}\label{th:delta2Isigman}
For positive integers $n$, the following are equivalent for quantifier-free $A,B$ and a free variable $c$.
\begin{enumerate}
\item
$I\Sigma^{0}_{n}$ proves 
\begin{equation}\label{eq:delta2}
\forall x\exists y A(x,y,c) \leftrightarrow \exists z\forall u B(z,u,c)
\end{equation}

\item
There exists a binary {\rm elementary} recursive predicate $f$, an ordinal $K<\omega_{1+n}$ and an {\rm elementary} recursive function
$h:\omega \times\omega\to K$  such that

\begin{enumerate}
\item  {\rm (weakly descending)}
\[
\mbox{{\rm {\sf EA}}}\vdash K>h(c,w)\geq h(c,w+1)
\]

\item {\rm (lowering)}
\[
\mbox{{\rm {\sf EA}}}\vdash f(c,w)\neq f(c,w+1) \to h(c,w)>h(c,w+1)
\]

\item 
\begin{eqnarray*}
\mbox{{\sf EA}} & \vdash & \lim_{w\to\infty}f(c,w)=0 \to \exists z\forall u B(z,u,c)
\\
\mbox{{\sf EA}} & \vdash & \lim_{w\to\infty}f(c,w)=1 \to \exists x\forall y \lnot A(x,y,c)
\\
I\Sigma^{0}_{n} & \vdash & \exists z\forall u B(z,u,c) \to \forall x\exists y A(x,y,c)
\end{eqnarray*}

\end{enumerate}
where the ordering $<$ denotes a standard $\varepsilon_{0}$-ordering $<_{\varepsilon_{0}}$.
\end{enumerate}
\end{theorem}

First note that by Theorem \ref{th:Sommer} we have $\Sigma^{0}_{1}$-minimization up to each
ordinal less than $\omega_{1+n}$
in $I\Sigma^{0}_{n}$. 
Hence $\exists\alpha<_{\varepsilon_{0}}K[\alpha=\min_{<_{\varepsilon_{0}}}\{\beta: \exists w[\beta=h(c,w)]\}]$.
Pick a $w$ so that the least $\alpha=h(c,w)$.
Assuming that {\rm {\sf EA}} (a fortiori $I\Sigma^{0}_{n}$) 
proves (weakly descending) and (lowering),
we have
\[
I\Sigma^{0}_{n} \vdash \forall u\geq w[f(c,u)=f(c,w)]
.\]
Therefore the convergence of the predicate $f$ is shown in $I\Sigma^{0}_{n}$.
Also
\[
I\Sigma^{0}_{n} \vdash \forall x\exists y A(x,y,c)\leftrightarrow \lim_{w\to \infty}f(c,w)=0
.\]

The converse follows from the following Reduction Theorem \ref{th:reduction}.

The theorem says that if a disjunction $\exists x\forall y \lnot A(x,y,c)\lor \exists z\forall u B(z,u,c)$ 
of $\Sigma^{0}_{2}$-formulas is provable in $I\Sigma^{0}_{n}$,
then one can construct an elementary recursive predicate $f$ whose limit tells us which disjunct is true.
The convergence of $f$ is ensured by a descending function $h$ in ordinals$<\omega_{1+n}$.
Moreover these are all provable in {\sf EA}.

Assuming the convergence of $f$(, which is provable in $I\Sigma^{0}_{n}$)
this is a provable version of the classical Reduction Property of $\Sigma^{0}_{2}$ sets to $\Delta^{0}_{2}$ sets.
The point is that the $\Delta^{0}_{2}$ sets $\{c:\lim_{w\to\infty}f(c,w)=0\}$ 
are in a level $D_{<\omega_{1+n}}$ of Ershov hierarchy, demonstrably in $I\Sigma^{0}_{n}$.

\begin{theorem}\label{th:reduction}{\rm (Reduction Property)}
Let $n\geq 1$.

Suppose $I\Sigma^{0}_{n}\vdash\exists x\forall y \lnot A(x,y,c)\lor \exists z\forall u B(z,u,c)$
for quantifier-free $A,B$.
Then there exists an {\rm elementary} recursive predicate $f$, an ordinal $K<\omega_{1+n}$ and an {\rm elementary} recursive function $h$ such that
\begin{enumerate}
\item  {\rm (weakly descending)}
\[
\mbox{{\rm {\sf EA}}}\vdash K>h(c,w)\geq h(c,w+1)
\]

\item {\rm (lowering)}
\[
\mbox{{\rm {\sf EA}}}\vdash f(c,w)\neq f(c,w+1) \to h(c,w)>h(c,w+1)
\]

\item {\rm (reduction)}
\begin{eqnarray*}
\mbox{{\sf EA}} & \vdash & \lim_{w\to\infty}f(c,w)=0 \to \exists z\forall u B(z,u,c)
\\
\mbox{{\sf EA}} & \vdash & \lim_{w\to\infty}f(c,w)=1 \to \exists x\forall y \lnot A(x,y,c)
\end{eqnarray*}
\end{enumerate}
\end{theorem}

In what follows, given a $I\Sigma^{0}_{n}$-proof of $\exists x\forall y \lnot A(x,y,c)\lor \exists z\forall u B(z,u,c)$
 let us construct a predicate $f$, an ordinal $K<\omega_{1+n}$ and a function $h$ 
enjoying (weakly descending), (lowering) and (reduction).

Let $p(x,y,c)$ denote the characteristic function of the predicate 
\[
A((x)_{0},(y)_{0},c)\to B((x)_{1},(y)_{1},c)
,\]
where $(x)_{i}\, (i=0,1)$ denotes the projections of the pairing function.

Then
\[
\exists x\forall y[ p(x,y,c)=0]
\]
is provable in $I\Sigma^{0}_{n}$.

\subsection{Infinitary derivations}

In what follows let us consider (finite or infinite) derivations in one-sided sequent calculi.
Given a finite derivation of $\exists x\forall y[ p(x,y,c)=0]$ in $I\Sigma^{0}_{n}$, first eliminate cut inferences partially to get a derivation of the same formula in which any cut formula is $\Sigma^{0}_{n}$.

Next embed the derivation into an infinite derivation of the sentence
\[
\exists x\forall y[ p(x,y,\bar{c})=0]
\]
 with the $c$-th numeral $\bar{c}$.
Then eliminate cut inferences to get a cut-free derivation $P_{c}$ of the same sentence.
As usual the depth of $P_{c}$ is bounded by an ordinal $K<\omega_{1+n}$ uniformly,
i.e., $\forall c[\mbox{dp}(P_{c})<K]$.

In the derivation $P_{c}$, the {\it initial sequents\/} are
\[
(Int)\:\: \Gamma, E
\]
for true equation $E$.
The equation $E$ is called the {\it main formula\/} of the initial sequent.

In what follows we identify the closed term $t$ with the numeral $\bar{n}$ of its value $n=val(t)$.

Note that 
the value of closed terms and truth values of equations in {\sf LEA} are elementary
recursively computable.
The initial sequents are regarded as inference rules with empty premiss (upper sequent), and with the empty list of side formulas.

The {\it inference rules\/} are $(\exists), (\forall)$, and the repetition rule $(Rep)$.
These are standard ones.
\[
\infer[(\exists)]
{\Gamma, \exists x B(x)}
{\Gamma, B(\bar{n})}
\: ;\:
 \infer[(\forall)]
 {\Gamma, \forall x B(x)}
 {
 \cdots
 &
 \Gamma,B(\bar{n})
 &
 \cdots (n\in\omega)
 }
\: ;\:
 \infer[(Rep)]
 {\Gamma}{\Gamma}
 \]
where $\exists x B(x)$ in the $(\exists)$ and $\forall x B(x)$ in the $(\forall)$ are the {\it main formula\/} of the inference,
and $B(\bar{n})$ are {\it side formulas\/} of the inferences.
The inference $(Rep)$ has no main nor side formulas.

Our infinitary derivations are equipped with additional informations as in \cite{Mintsfinite}.
\begin{definition}\label{df:derivation}
{\rm An} infinitary derivation {\rm is a sextuple} 
\[
D=(T,Seq,Rule,Mfml,Sfml,ord)
\]
 {\rm which enjoys the following conditions.
The naked tree of} $D$ {\rm is denoted} $T=T(D)$.
\begin{enumerate}
\item
$T\subseteq{}^{<\omega}\omega$ {\rm is a tree with its root} $\emptyset$ {\rm such that}
\[
a*\langle n\rangle\in T \,\&\, m<n \Rightarrow a*\langle m\rangle\in T
.\]

\item
$Seq(a)$ {\rm for} $a\in T$ {\rm denotes the sequent situated at the node} $a$.

{\rm If} $Seq(a)$ {\rm is a sequent} $\Gamma${\rm , then it is denoted}
\[
a:\Gamma
.\]

\item
$Rule(a)$ {\rm for} $a\in T$ {\rm denotes the name of the inference rule with its lower sequent} $Seq(a)$.

\item
$Mfml(a)$  {\rm for} $a\in T$ {\rm denotes the} main formula {\rm of the inference rule} $Rule(a)$.
{\rm When} $Rule(a)=(Rep)${\rm , then} $Mfml(a)=\emptyset$.

\item
$Sfml(a*\langle n\rangle)$  {\rm for} $a*\langle n\rangle\in T$ {\rm denotes the} side formula {\rm of the inference rule} $Rule(a)${\rm , which is in the} $n${\rm -th upper sequent, i.e.,} $Sfml(a*\langle n\rangle)\in Seq(a*\langle n\rangle)$.
{\rm When} $Rule(a)=(Rep), (Int)${\rm , then} $Sfml(a*\langle n\rangle)=\emptyset$.

\item
$ord(a)$ {\rm for} $a\in T$ {\rm denotes the ordinal}$<_{\varepsilon_{0}}K$ {\rm attached to} $a$.

\item
{\rm The sextuple} $(T,Seq,Rule,Mfml,Sfml,ord)$ {\rm has to be locally correct with respect to
inference rules of the infinitary calculus and
for being well founded tree} $T$.

\end{enumerate}
\end{definition}

In a derivation each inference rule except $(Int)$ receives the following nodes:
\[
\infer[(\exists)]
{a: \Gamma, \exists x B(x)}
{a*\langle 0\rangle : \Gamma, B(\bar{n})}
\: ;\:
 \infer[(\forall)]
 {a: \Gamma, \forall x B(x)}
 {
 \cdots
 &
 a*\langle n\rangle : \Gamma,B(\bar{n})
 &
 \cdots (n\in\omega)
 }
\: ;\:
 \infer[(Rep)]
 {a: \Gamma}{a*\langle 0\rangle: \Gamma}
 \]

The ordinals $ord_{c}(a)$ in the inference $(\forall)$
\[
 \infer[(\forall)]
 {a: \Gamma, \forall x B(x)}
 {
 \cdots
 &
 a*\langle n\rangle : \Gamma,B(\bar{n})
 &
 \cdots (n\in\omega)
 }\]
enjoys
\begin{equation}\label{eq:ordass}
ord_{c}(a)>ord_{c}(a*\langle n\rangle)=ord_{c}(a*\langle m\rangle)
\end{equation}
for any $n,m$.

As in \cite{Mintsfinite} we see that the function $c\mapsto P_{c}$ is elementary recursive.
We denote $P_{c}=(T_{c},Seq_{c},Rule_{c},Mfml_{c},Sfml_{c},ord_{c})$.

\subsection{Searching witnesses of $\Sigma^{0}_{2}$ in derivations}

Let us define a tracing function $\sigma(c,i)\in T_{c}=T(P_{c})$.

The function $\{\sigma(c,w)\}_{w}$ indicates the trail in the proof tree $T_{c}$
in which we go through in searching a witness $x_{a}$ of $\exists x\forall y[ p(x,y,\bar{c})=0]$, and verifying 
$\forall y[p(x_{a},y,\bar{c})=0]$.

\begin{enumerate}
\item
$\sigma(c,0)=\emptyset$(root).

In what follows let $a=\sigma(c,w)$.

\item
Until $Seq_{c}(a)$ is an upper sequent of an $(\forall)$, go to the leftmost branch:
\[
\sigma(c,w+1)=a*\langle 0\rangle
.\]
For example
\[
\infer[(\exists)]
{a: \Gamma,\exists x\forall y[p(x,y,\bar{c})=0]}
{
a*\langle 0\rangle : \Gamma,\exists x\forall y[p(x,y,\bar{c})=0],\forall y[p(x_{a},y,\bar{c})=0]
}
\]

\item
The case when $Rule_{c}(b)=(\forall)$ with $a=b*\langle n\rangle$. 
Namely $Seq_{c}(a)$ is the $n$-th upper sequent of an $(\forall)$.
\[
\infer[(\forall)]
{\Gamma,\forall y[p(x_{a},y,\bar{c})=0]}
{
\cdots
&
a: \Gamma,p(x_{a},y_{a},\bar{c})=0
&
\cdots
}
\]
$x_{a},y_{a}$ are closed terms.

 \begin{enumerate}
 \item
  If $p(x_{a},y_{a},\bar{c})=0$ is a TRUE equation, 
  $\sigma(c,w+1)=a\oplus 1$, the next right to the $a$:
 \[
\infer[(\forall)]
{\Gamma,\forall y[p(x_{a},y,\bar{c})=0]}
{
\sigma(c,w) : \Gamma,p(x_{a},y_{a},\bar{c})=0
&
\sigma(c,w+1) : \Gamma,p(x_{a},y_{a}+1,\bar{c})=0
}
\] 
where for an $a=(a_{0},\ldots,a_{n-2},a_{n-1})\in{}^{<\omega}\omega$
  \[
  a\oplus 1=(a_{0},\ldots,a_{n-2},a_{n-1}+1)
  \]
  if $lh(a)=n>0$.
  
  $\emptyset\oplus 1$ is defined to be $\emptyset$.

 \item
  Otherwise  $\sigma(c,w+1)=a*\langle 0\rangle$, i.e., go to the leftmost branch from $a$.
 
\[
\infer[(\forall)]
{\Gamma,\forall y[p(x_{a},y,\bar{c})=0]}
{
\cdots
&
\infer{\sigma(c,w) : \Gamma,p(x_{a},y_{a},\bar{c})=0}{\sigma(c,w+1):\Delta & \cdots}
&
\cdots
}
\]

 \end{enumerate}

It is easy to see that the function $(c,w)\mapsto \sigma(c,w)$ is elementary recursive
since $\max(\{(\sigma(c,w))_{i}: i<lh(\sigma(c,w))\}\cup\{lh(\sigma(c,w))\})\leq w$.
\end{enumerate}

Once $\sigma(c,w)$ is on an $(\forall)$, the tracing function goes through the upper sequents as long as the equations $p(x_{a},y_{a},\bar{c})=0$ is TRUE.

It is intuitively clear that after a finite number of steps, 
the sequence $\{\sigma(c,w)\}_{w}$ goes through the upper sequents of an $(\forall)$:

 \[
\infer[(\forall)]
{\Gamma,\forall y[p(x_{a},y,\bar{c})=0]}
{
\sigma(c,w_{0}) : \Gamma,p(x_{a},0,\bar{c})=0
&
\cdots 
&
\sigma(c,w_{0}+y) : \Gamma,p(x_{a},\bar{y},\bar{c})=0
&
\cdots
}
\] 
since $\forall y[p(x_{a},y,\bar{c})=0]$ is true for an $x_{a}$.
We will know at the limit the fact, i.e.,
for $x=(x_{a})_{0}$ and $z=(x_{a})_{1}$
\[
\exists y A(\bar{x},y,\bar{c})\to \forall u B(\bar{z},u,\bar{c})
\]
is true.

Now let us define an elementary recursive predicate $f$ as follows.

\begin{enumerate}
\item
$f(c,0)=1$.

\item
Alternate values
$f(c,w+1)=1-f(c,w)$ if
$Seq_{c}(\sigma(c,w+1))$ is an upper sequent of an inference other than $(\forall)$.

\item
Suppose $Seq_{c}(\sigma(c,w+1))$ is the $n$-th upper sequent of an $(\forall)$, and
$\sigma(c,w+1)=b*\langle n\rangle$.

 \[
\infer[(\forall)]
{b: \Gamma,\forall y,u[A(x_{b},y,\bar{c})\to B(z_{b},u,\bar{c})]}
{
\cdots 
&
b*\langle n\rangle : \Gamma, A(x_{b},(\bar{n})_{0},\bar{c})\to B(z_{b},(\bar{n})_{1},\bar{c})
&
\cdots
}
\] 

$f(c,w+1)=0$ iff $A(x_{b},(n)_{0},\bar{c})\to B(z_{b},(n)_{1},\bar{c})$ is true, and
the following condition holds:
\[
  \exists k\leq n [A(x_{b}, (\bar{k})_{0},\bar{c})]
  \]

\end{enumerate}

Namely
\begin{eqnarray*}
&&
f(c,w+1)=0
\Leftrightarrow
\\
&&
[A(x_{b},(n)_{0},\bar{c})\to B(z_{b},(n)_{1},\bar{c})] \,\&\,
  \exists k\leq n [A(x_{b}, (\bar{k})_{0},\bar{c})]
\end{eqnarray*}
  
Suppose $\sigma(c,w)$ is on an $(\forall)$.
Until a witness $k$ such that $A(x_{b}, (\bar{k})_{0},\bar{c})$ is found, $f(c,w)=1\, (w<k)$.
After a witness $k$ has been found, $f(c,w)=0\, (w\geq k)$ as long as 
 $A(x_{b},(\bar{n})_{0},\bar{c})\to B(z_{b},(\bar{n})_{1},\bar{c})$ is true.

Therefore if the tracing function $\sigma(c,w)$ goes through the upper sequents of the $(\forall)$, then
either $\lim_{w\to\infty}f(c,w)=1$ and $\forall y \lnot A(x_{b}, y,\bar{c})$, 
or $\lim_{w\to\infty}f(c,w)=0$ and $\forall u B(z_{b},u,\bar{c})$.

\begin{proposition}\label{prp:delta2Isigman.1}
\begin{enumerate}
\item\label{prp:delta2Isigman.11}
Suppose that $b*\langle n\rangle=\sigma(c,w+1)$ and $Seq_{c}(b*\langle n\rangle)$ is the $n$-th upper sequent of an inference
$(\forall)$.
Then $\{f(c,u): \sigma(c,u)=b*\langle m\rangle, m\leq n\}$ changes the values at most twice.
Moreover if $f(c,u)=0$ and $f(c,v)=1$ for some $u<v\leq w+1$, then
$v=w+1$ and $\sigma(c,v+1)=\sigma(c,v)*\langle 0\rangle$, i.e., $Seq_{c}(\sigma(c,v))$ is the last upper sequent
of the inference $(\forall)$ in the tracing function $\sigma$.
\item\label{prp:delta2Isigman.12}{\rm (Reduction)}
\begin{eqnarray*}
\mbox{{\sf EA}} & \vdash & \lim_{w\to\infty}f(c,w)=0 \to \exists z\forall u B(z,u,c)
\\
\mbox{{\sf EA}} & \vdash & \lim_{w\to\infty}f(c,w)=1 \to \exists x\forall y \lnot A(x,y,c)
\end{eqnarray*}
\end{enumerate}
\end{proposition}
{\bf Proof}.\hspace{2mm}
Recall that an inference rule $(\forall)$ in $P_{c}$ is of the form:
\[
\infer[(\forall)]
{b: \Gamma, \forall y [p(x_{b}, y,\bar{c})=0]}
{
\ldots
&
b*\langle n\rangle: \Gamma,p(x_{b*\langle n\rangle},\bar{n},\bar{c})=0
&
\ldots
}
\]
where
\[
p(x_{b*\langle n\rangle},\bar{n},\bar{c})=0 \leftrightarrow
[A(((x)_{b*\langle n\rangle})_{0}, (\bar{n})_{0},\bar{c}) \to 
B((x)_{b*\langle k\rangle})_{1}, (\bar{n})_{1},\bar{c}) ]
\]

Let $u$ be such that $\sigma(c,u)=b*\langle m\rangle$ with an $m\leq n$.
Then by the definition of the tracing function $\sigma$, we have for $m<n$
$p(x_{b*\langle m\rangle},\bar{m},\bar{c})=0$, i.e.,
\[
A((x_{b*\langle m\rangle})_{0},(\bar{m})_{0},\bar{c})\to B((x_{b*\langle m\rangle})_{1},(\bar{m})_{1},\bar{c})
.\]
Suppose there exists a $u\leq w+1$ such that $f(c,u)=0$, and let $u$ denote the minimal such one.

Then for any $v$ with $u\leq v<w+1$, we have $f(c,v)=0$.
Therefore if $f(c,v)=1$ for a $v>u$, it must be the case $v=w+1$.
This means that for some  $k\leq n-1$ with $(k)_{0}=(n)_{0}$
\[
A((x_{b})_{0},(k)_{0},c)\land \lnot B((x_{b})_{1},(n)_{1},c)
.\]
Hence $p(x_{b},n,c)\neq 0$, and $\sigma(c,v+1)=\sigma(c,v)*\langle 0\rangle$.
\hspace*{\fill} $\Box$

Next define $h$ as follows.
\begin{enumerate}
\item
\[
h(c,0)=3\cdot ord_{c}(\emptyset)
.\]

In what follows put $a=\sigma(c,w+1)$ and let $Seq_{c}(a)$ be an upper sequent of an inference
$Rule_{c}(b)$ with $a=b*\langle n\rangle$.

\item
The case when $Rule_{c}(b)$ is an inference rule other than $(\forall)$.
\[
h(c,w+1):=3\cdot ord_{c}(\sigma(c,w+1))
.\]

By Proposition \ref{prp:delta2Isigman.1}.\ref{prp:delta2Isigman.11} we know that
the $f(c,u)$ changes the values at most twice in the upper sequents of an $(\forall)$.
\item
The case when $n=0$ and $Rule_{c}(b)=(\forall)$.
\[
h(c,w+1):=3\cdot ord_{c}(\sigma(c,w+1))+2
.\]

\item
The case when $n>0$, $Rule_{c}(b)=(\forall)$.
 \[
\infer[(\forall)]
{b: \Gamma,\forall y,u[A(x_{b},y,\bar{c})\to B(z_{b},u,\bar{c})]}
{
\cdots 
&
b*\langle n\rangle : \Gamma, A(x_{b},(\bar{n})_{0},\bar{c})\to B(z_{b},(\bar{n})_{1},\bar{c})
&
\cdots
}
\] 
We have by (\ref{eq:ordass})
\[
ord_{c}(\sigma(c,w))=ord_{c}(\sigma(c,w+1))
.\]

 \begin{enumerate}
  \item
The case when $f(c,w+1)=f(c,w)$.
\[
h(c,w+1):=h(c,w)
.\]
where $\sigma(c,w+1)=\sigma(c,w)\oplus 1$.

 \item
The case when $f(c,w)=1\,\&\, f(c,w+1)=0$.

Then
$\sigma(c,w+1)=\sigma(c,w)\oplus 1$ and 
$n=\min\{k: A(x_{b}, (\bar{k})_{0},\bar{c})\}$.

Let
\[
h(c,w+1):=3\cdot ord_{c}(\sigma(c,w+1))+1
.\]

 \item
 The case when $f(c,w)=0\,\&\, f(c,w+1)=1$.
 
 This means that $A(x_{b}, (\bar{n})_{0},\bar{c})\to B(z_{b}, (\bar{n})_{1},\bar{c})$ is FALSE
and
\\ $\sigma(c,w+2)=\sigma(c,w+1)*\langle 0\rangle$.
 \[
h(c,w+1):=3\cdot ord_{c}(\sigma(c,w+1))
.\]

 \end{enumerate}

Obviously $h$ is elementary recursive.

\end{enumerate}

\begin{proposition}\label{prp:delta2Isigman.2}

\[
\mbox{{\rm (weakly descending) }}
\mbox{{\sf EA}} \vdash 3K>h(c,w)\geq h(c,w+1)
\]

\[
\mbox{{\rm (lowering) }}
\mbox{{\sf EA}} \vdash f(c,w)\neq f(c,w+1) \to h(c,w)>h(c,w+1)
\]
\end{proposition}
{\bf Proof}.\hspace{2mm}
 (weakly descending) is obvious.
 
Consider the case when $\sigma(c,w+1)=a$ and $Seq_{c}(a)$ is an upper sequent of an inference
$Rule_{c}(b)=(\forall)$ with $a=b*\langle n\rangle$.

If $n=0$, then 
\[
h(c,w+1)=3\cdot ord_{c}(\sigma(c,w+1))+2<3\cdot ord_{c}(\sigma(c,w))\leq h(c,w)
\]
since $Seq_{c}(\sigma(c,w))$ is the lower sequent of $Seq_{c}(a)$ with
$b=\sigma(c,w)$.

Assume $n>0$.
Using Proposition \ref{prp:delta2Isigman.1}.\ref{prp:delta2Isigman.11} we see
$h(c,w+1)\in\{3\cdot ord_{c}(\sigma(c,w+1))+i: i<3\}$.

Moreover if $\sigma(c,w+2)=\sigma(c,w+1)*\langle 0\rangle$, then
\[
h(c,w+1)\geq 3\cdot ord_{c}(\sigma(c,w+1))>3\cdot ord_{c}(\sigma(c,w+1))+2\geq h(c,w+2)
.\]

\hspace*{\fill} $\Box$

This completes a proof of Theorems \ref{th:reduction} and \ref{th:delta2Isigman}.

\subsection{Provably $\Sigma^{0}_{2}$-functions}

If $\exists z\forall u B(z,u,c)$ is provable for quantifier-free $B$, 
then we can find a witness $z= \lim_{w\to\infty}f(c,w)$ as limits of an $f$,
whose convergence is ensured by weakly descending chains of ordinals.

\begin{theorem}\label{th:sigma2witness}

Suppose $I\Sigma^{0}_{n}\vdash \exists z\forall u B(z,u,c)$
for quantifier-free $B$.
Then there exist {\rm elementary} recursive functions $f$, $h$ and an ordinal $K<\omega_{1+n}$ such that
\begin{enumerate}
\item 
\[
\mbox{{\rm (weakly descending) }}
\mbox{{\sf EA}} \vdash K>h(c,w)\geq h(c,w+1)
\]

\item
\[
\mbox{{\rm (lowering) }}
\mbox{{\sf EA}} \vdash f(c,w)\neq f(c,w+1) \to h(c,w)>h(c,w+1)
\]

\item 
\[
\mbox{{\sf EA}}  \vdash  \lim_{w\to\infty}f(c,w)=z \to \forall u B(z,u,c)
\]

\end{enumerate}

\end{theorem}
{\bf Proof}.\hspace{2mm}
As in the proof of Theorem \ref{th:reduction}, let us define a tracing function $\sigma$.

$\sigma(c,w)$ goes on the leftmost branch up to an $(\forall)$.
$\sigma(c,w)$ goes through the upper sequents of $(\forall)$ as long as side formulas $B(z_{a},\bar{n},\bar{c})$ is TRUE.
If a FALSE side formula $B(z_{a},\bar{n},\bar{c})$ is found, then throw $z_{a}$ away and go on the leftmost branch.

Now $h$ is defined by $h(c,w):=ord_{c}(\sigma(c,w))$.
 $f$ is defined obviously.
 $f(c,w)=z_{a}$ if $Seq_{c}(\sigma(c,w))$ is an upper sequent of an $(\forall)$ with its side formula $B(z_{a},\bar{n},\bar{c})$.
 Otherwise $f(c,w)$ is arbitrary, say $f(c,w)=0$.
 \hspace*{\fill} $\Box$
 \\ \smallskip \\
It is well known that the 1-consistency $\mbox{RFN}_{\Pi^{0}_{2}}(I\Sigma^{0}_{n})$
is equivalent over Primitive Recursive Arithmetic {\sf PRA} to the fact that there is no primitive recursive descending chain
of ordinals$<\omega_{1+n}$.

\begin{theorem}\label{th:2con}(Cf. \cite{AraiMints} for another form of the 2-consistency of arithmetic.)

The 2-consistency $\mbox{{\rm RFN}}_{\Pi^{0}_{3}}(I\Sigma^{0}_{n})$
is equivalent over {\sf PRA} to the fact that every primitive recursive 
weakly descending chain of ordinals$<\omega_{1+n}$ has a limit, 
or equivalently to the fact that
for any primitive recursive sequence $\{h(c,w)\}_{w}$ of ordinals$<\omega_{1+n}$ 
the least ordinal $\min_{<_{\varepsilon_{0}}}\{h(c,w): w\in\omega\}$ exists.
\end{theorem}
{\bf Proof}.\hspace{2mm}

Over {\sf PRA}, $\mbox{RFN}_{\Pi^{0}_{3}}(I\Sigma^{0}_{n})$ yields the existence of 
the least ordinal $\min_{<_{\varepsilon_{0}}}\{h(c,w)<\omega_{1+n}: w\in\omega\}$ since $\alpha=\min_{<_{\varepsilon_{0}}}\{\beta: \exists w[\beta=h(c,w)]\}$ is a $\Sigma^{0}_{2}$-formula.

Conversely let $f(c,w)<2$ be defined as follows:
\begin{enumerate}
\item
$c$ is not a G\"odel number of an $I\Sigma^{0}_{n}$-proof of a $\Sigma^{0}_{2}$-sentence:
Then $f(c,w)=0$ for any $w$.

\item
$c$ is a G\"odel number of an $I\Sigma^{0}_{n}$-proof of a $\Sigma^{0}_{2}$-sentence $\exists z\forall u B_{c}(z,u)$:
\\
$f(c,w)$ is defined as in Theorem \ref{th:sigma2witness} for
a cut free infinite derivation $P_{c}$ of $\exists z\forall u B_{c}(z,u)$.
Note that $f$ is non-elementary since it involves cut elimination for predicate logic.
\end{enumerate}
Also let
 $h(c,w):=ord_{c}(\sigma(c,w))$.
 
Then

\begin{enumerate}
\item 
\[
\mbox{{\rm (weakly descending) }}
\mbox{{\sf PRA}} \vdash \omega_{1+n}>h(c,w)\geq h(c,w+1)
\]

\item
\[
\mbox{{\rm (lowering) }}
\mbox{{\sf PRA}} \vdash f(c,w)\neq f(c,w+1) \to h(c,w)>h(c,w+1)
\]

\item 
\[
\mbox{{\sf PRA}}  \vdash  \lim_{w\to\infty}f(c,w)=0 \to \mbox{{\rm Prov}}_{I\Sigma^{0}_{n}}(c,\lceil \exists z\forall u B_{c}(z,u)\rceil)\to \exists z\forall u B_{c}(z,u)
\]
\end{enumerate}
Therefore
\[
\mbox{{\sf PRA}}\vdash \forall c [\exists w\forall u\geq w \{h(c,u)=h(c,w)\} 
\to \exists\ell\{ \lim_{w\to\infty}f(c,w)=\ell\} ]
\]
and
\[
\mbox{{\sf PRA}}\vdash \forall c\exists\ell[ \lim_{w\to\infty}f(c,w)=\ell]\to \mbox{RFN}_{\Pi^{0}_{3}}(I\Sigma^{0}_{n})
.\]

\hspace*{\fill} $\Box$

\section{Provably $\Delta^{0}_{2}$ in {\sf EA}
}\label{sect:PRA}

In this section we consider the $\Delta^{0}_{2}$-sets provably in {\sf EA}.

The following Theorem \ref{th:delta2pra} states that for 
a set is provably $\Delta^{0}_{2}$ in {\sf EA}
 iff
it is {\sf EA}-provably in the class $D_{n}$ of a finite level in the Ershov hierarchy.
The finite levels $\{D_{n} : n<\omega\}$ are called the difference (or Boolean) hierarchy, and 
by a result due to H. Putnam(Theorem 2 in \cite{Putnam}) we see that
a set is provably $\Delta^{0}_{2}$ in {\sf EA} iff it is equivalent to a Boolean combination of $\Sigma^{0}_{1}$-formulas, provably in {\sf EA}.
This answers to a problem of L. Beklemishev.

\begin{theorem}\label{th:delta2pra}
The following are equivalent for quantifier-free $A,B$ and a free variable $c$.
\begin{enumerate}
\item
{\rm {\sf EA}} proves
\begin{equation}
\renewcommand{\theequation}{\ref{eq:delta2}}
\forall x\exists y A(x,y,c) \leftrightarrow \exists z\forall u B(z,u,c)
\end{equation}
\addtocounter{equation}{-1}

\item
There exists a binary {\rm elementary} recursive predicate $f$, a natural number $K<\omega$ and an {\rm elementary} recursive function
$h:\omega \times\omega\to K$  such that

\begin{enumerate}
\item  {\rm (weakly descending)}
\[
\mbox{{\rm {\sf EA}}}\vdash K>h(c,w)\geq h(c,w+1)
\]

\item {\rm (lowering)}
\[
\mbox{{\rm {\sf EA}}}\vdash f(c,w)\neq f(c,w+1) \to h(c,w)>h(c,w+1)
\]

\item {\rm (reduction)}
\begin{eqnarray*}
\mbox{{\sf EA}} & \vdash & \lim_{w\to\infty}f(c,w)=0 \to \exists z\forall u B(z,u,c)
\\
\mbox{{\sf EA}} & \vdash & \lim_{w\to\infty}f(c,w)=1 \to \exists x\forall y \lnot A(x,y,c)
\\
\mbox{{\sf EA}} & \vdash & \exists z\forall u B(z,u,c) \to \forall x\exists y A(x,y,c)
\end{eqnarray*}

\end{enumerate}
for the usual ordering $<$ on $\omega$.
\end{enumerate}
\end{theorem}
{\bf Proof}.\hspace{2mm}
Assume {\sf EA} proves (weakly descending) and (lowering) for a natural number $K$. Then {\sf EA} also proves the convergence of $f$:
\[
\mbox{{\rm {\sf EA}}}\vdash\exists \ell[\lim_{w\to\infty}f(c,w)=\ell]
,\]
(reduction) yields
\[
\mbox{{\rm {\sf EA}}}\vdash \forall x\exists y A(x,y,c)\leftrightarrow \lim_{w\to \infty}f(c,w)=0
.\]

Conversely
suppose that {\sf EA} proves (\ref{eq:delta2}).
Then so is the $\exists\forall$-formula
\[
\exists x\exists z\forall y\forall z[A(x,y,c)\to B(z,u,c)]
.\]
By the Herbrand's theorem there exist a list of variables $\{a_{i},b_{i}: i\leq r\}$ and
a list of terms $\{t_{i},s_{i}:i\leq r\}$ such that
\begin{equation}\label{eq:Herbrand}
\bigvee\{ A(t_{i}, a_{i},c) \to B(s_{i},b_{i},c) :i\leq r\}
\end{equation}
is provable in {\sf EA}, and variables occurring in $t_{i},s_{i}$ are among $a_{j},b_{j}$ for  $j<i$
besides the parameter $c$.

For simplicity consider the case when $r=1$. Then we have

\begin{equation}\label{eq:Herbrand1}
\mbox{{\sf EA}}\vdash \lnot A(t_{0}, a_{0},c) \lor B(s_{0},b_{0},c) \lor 
\lnot A(t_{1}(a_{0},b_{0}), a_{1},c) \lor B(s_{1}(a_{0},b_{0}),b_{1},c)
\end{equation}

Let $f$ denote the elementary recursive predicate
\[
f(c,w):=\left\{
\begin{array}{ll}
0 &  [\{t_{0}\leq w \land \exists y\leq w A(t_{0},y,c)\} \land
\{ s_{0}\leq w \land \forall u\leq w B(s_{0},u,c)\}] \lor 
\\
&
[\exists a_{0},b_{0}\leq w\{A(t_{0},a_{0},c)\land \lnot B(s_{0},b_{0},c) \land 
\\
&
\exists a_{1}\leq w  A(t_{1}(a_{0},b_{0}), a_{1},c) \land \forall b_{1}\leq w B(s_{1}(a_{0},b_{0}),b_{1},c)\}]
\\
1 & \mbox{ otherwise}
\end{array}
\right.
\]

For the number
\[
K:=1+2r+2 (=5 \mbox{ if } r=1)
,\]
let $h:\omega\times\omega\to K$ denote the elementary recursive function
\[
h(c,0)=K-1
\]
and
\[
h(c,w+1):=\left\{
\begin{array}{ll}
h(c,w) & \mbox{ if }  f(c,w+1)=f(c,w)
\\
h(c,w)\dot{-}1 & \mbox{ if } f(c,w+1)\neq f(c,w)
\end{array}
\right.
\]

\begin{lemma}\label{lem:1}
{\rm {\sf EA}} proves the facts 
(weakly descending), (lowering) and (reduction).
\end{lemma}
{\bf Proof}.\hspace{2mm}
Argue in {\sf EA}.
(weakly descending) is obvious.

Suppose
\[
\lim_{w\to\infty}f(c,w)=\ell
\]
for an $\ell=0,1$.

By (\ref{eq:delta2})
we have either $\exists z\forall u B(z,u,c)$ or $\exists x\forall y \lnot A(x,y,c)$.

First consider the case when $\exists z\forall u B(z,u,c)$. Then $\forall x\exists y A(x,y,c)$.
Hence by (\ref{eq:Herbrand1}) either $\forall b_{0} B(s_{0},b_{0},c)$ or
$\forall b_{1} B(s_{1}(a_{0},b_{0}),b_{1},c)$ for some $a_{0}, b_{0}$
with $A(t_{0},a_{0},c)\land \lnot B(s_{0},b_{0},c)$.

If $\forall b_{0} B(s_{0},b_{0},c)$, then $f(c,w)=0$ for any $w\geq\max\{t_{0},s_{0},y_{0}\}$,
where $y_{0}=\mu y. A(t_{0},y.c)$.
Therefore $\ell=0$.
Moreover $f(c,w)=1$ for $w<\max\{t_{0},s_{0},y_{0}\}$.

Next assume $\forall b_{1} B(s_{1}(a_{0},b_{0}),b_{1},c)$ for the minimal $a_{0}.b_{0}$
such that $A(t_{0},a_{0},c)\land \lnot B(s_{0},b_{0},c)$.
Then let $a_{1}$ denote the minimal $a_{1}$ such that $A(t_{1}(a_{0},b_{0}),a_{1},c)$.
We have $f(c,w)=0$ for any $w\geq\max\{a_{0},b_{0},a_{1}\}$, and hence $\ell=0$.

Now consider $w<\max\{a_{0},b_{0},a_{1}\}$.
Then $f(c,w)=0$ iff $\max\{t_{0},s_{0},a_{0}\}\leq w<b_{0}$.
Therefore $\lambda w.f(c,w)$ changes its values at most three times
(when $\max\{t_{0},s_{0},a_{0}\}<b_{0}<a_{1}$).

Next consider the case when $\exists x\forall y \lnot A(x,y,c)$.
We have $\forall z\exists u\lnot B(z,u,c)$.
Then $f(c,w)=1$ for any $w$
if $\forall a_{0}\lnot A(t_{0},a_{0},c)$,
and $f(c,w)=1$ for any $w\geq\max\{b_{0},b_{1}\}$
if $\forall a_{1}\lnot A(t_{1}(a_{0},b_{0}),a_{1},c)$
for the minimal $a_{0}, b_{0},b_{1}$
such that $A(t_{0},a_{0},c)\land \lnot B(s_{0},b_{0},c)$
and $\lnot B(s_{1}(a_{0},b_{0}),b_{1},c)$.
Therefore $\ell=1$.

Finally assume $\forall a_{1}\lnot A(t_{1}(a_{0},b_{0}),a_{1},c)$, and consider $w<\max\{b_{0},b_{1}\}$.
Then $f(c,w)=0$ iff $\max\{t_{0},s_{0},a_{0}\}\leq w<b_{0}$.
Therefore $\lambda w.f(c,w)$ changes its values at most two times in this case.

In any cases, (reduction) was shown, and $\lambda w.f(c,w)$ changes its values at most $1+2r(=3 \mbox{ {\rm if }} r=1)$ times for any $c$, i.e.,
\[
\forall (w_{0}<w_{1}<\cdots<w_{1+2r})\exists i\leq 1+2r [ f(c, w_{i})=f(c, w_{i}+1)]
.\]
Hence (lowering) follows.
\hspace*{\fill} $\Box$

Lemma \ref{lem:1} with a result due to H. Putnam(Theorem 2 in \cite{Putnam}) yields the
\begin{theorem}\label{th:1}
Suppose that {\sf EA} proves
\[
\forall x\exists y A(x,y,c) \leftrightarrow \exists z\forall u B(z,u,c)
\]
for quantifier-free $A,B$.

Then over {\sf EA} $\exists z\forall u B(z,u,c)$ is equivalent to a Boolean combination of $\Sigma^{0}_{1}$-formulas.
\end{theorem}
{\bf Proof}.\hspace{2mm}(cf. \cite{Putnam}.)
Let $r$ be as in (\ref{eq:Herbrand}).
Define $\Sigma^{0}_{1}$ $Y_{k}(c)\, (k\leq 1+2r)$ and $N_{i}(c)$ by
\begin{eqnarray*}
Y_{k}(c) & :\Leftrightarrow & \exists (w_{0}<w_{1}<\cdots<w_{k-1})\forall i<k [ f(c, w_{i})\neq f(c, w_{i}+1) \\
&& \land  f(c,w_{k-1})=0]
\\
N_{k}(c) & :\Leftrightarrow &
\exists (w_{0}<w_{1}<\cdots<w_{k-1})\forall i<k [ f(c, w_{i})\neq f(c, w_{i}+1)  \\
&& \land
f(c,w_{k-1})=1]]
\end{eqnarray*}
for $k>0$, and
$Y_{0}(c):\Leftrightarrow f(c,0)=0$, $N_{0}(c):\Leftrightarrow f(c,0)=1$.
Also put $N_{2r+2}(c):\Leftrightarrow 0=1$.

Then {\sf EA} proves that
\[
\exists z\forall u B(z,u,c) \leftrightarrow \bigvee\{Y_{k}(c)\land\lnot N_{k+1}(c): k\leq 1+2r\}
.\]
\hspace*{\fill} $\Box$

As in Theorems \ref{th:sigma2witness}, \ref{th:2con} we see the following theorems.

\begin{theorem}\label{th:sigma2witnesspra}

Suppose ${\sf EA}\vdash \exists z\forall u B(z,u,c)$
for quantifier-free $B$.
Then there exist {\rm elementary} recursive functions $f$, $h$ and a natural number $K<\omega$ such that
\begin{enumerate}
\item 
\[
\mbox{{\rm (weakly descending) }}
\mbox{{\sf EA}} \vdash K>h(c,w)\geq h(c,w+1)
\]

\item
\[
\mbox{{\rm (lowering) }}
\mbox{{\sf EA}} \vdash f(c,w)\neq f(c,w+1) \to h(c,w)>h(c,w+1)
\]

\item 
\[
\mbox{{\sf EA}}  \vdash  \lim_{w\to\infty}f(c,w)=z \to \forall u B(z,u,c)
\]

\end{enumerate}

\end{theorem}

\begin{theorem}\label{th:2conpra}
The 2-consistency $\mbox{{\rm RFN}}_{\Pi^{0}_{3}}({\sf EA})$
is equivalent over {\sf PRA} to the fact that every primitive recursive 
weakly descending chain of natural number$<\omega$ has a limit, 
or equivalently to the fact that
for any primitive recursive sequence $\{h(c,w)\}_{w}$ of natural number$<\omega$
the least number $\min_{<}\{h(c,w)<\omega: w\in\omega\}$ exists.
\end{theorem}
{\bf Remark}.

Obviously Theorems \ref{th:delta2pra}, \ref{th:1} and \ref{th:sigma2witnesspra} 
hold for any purely universal extension of {\sf EA}, eg., {\sf EA}+CON({\sf EA}), {\sf PRA}.

\appendix

\section{Nested limit existence rules}\label{appendix}

Every fragment in the Appendix is an extension of Elementary Recursive Arithmetic {\sf EA}.

In \cite {Lev}, Beklemishev and Visser gave an elegant axiomatization of $\Sigma^{0}_{2}$-consequences of
$I\Sigma^{0}_{1}$ in terms of the inference rule $(LimR)$ for limit existence principle:
\[
\infer[(LimR)]{\exists m\forall n\geq m \, h(n)=h(m)}{\exists m\forall n\geq m\, h(n+1)\leq h(n)}
\]

Moreover unnested applications of $(LimR)$ is shown to be equivalent to $I\Pi_{1}^{-}$ (over {\sf EA}).

This reminds us another axiomatization of $\Sigma^{0}_{2}$-consequences of
$I\Sigma^{0}_{1}$ in \cite{KPD}.
Namely $I\Sigma^{0}_{1}$ is a  $\Sigma^{0}_{2}$ conservative extension of 
$L\Sigma_{1}^{-(\infty)}=\bigcup_{k}L\Sigma_{1}^{-(k)}$, where
$L\Sigma_{1}^{-(k)}$ denotes the schema
\begin{eqnarray*}
&&
\exists x_{1}\cdots\exists x_{k} \theta(x_{1},\ldots,x_{k}) \to
\\
&&
\exists x_{1}\cdots\exists x_{k} \bigwedge_{i=1}^{k}[
\exists\vec{y} \theta(x_{1},\ldots,x_{i},\vec{y}) \land \forall z<x_{i}\forall\vec{y}\lnot\theta(x_{1},\ldots,x_{i-1},z,\vec{y})]
\end{eqnarray*}
for $\theta\in\Sigma^{0}_{1}$ without parameters.

For example $L\Sigma_{1}^{-(0)}={\sf EA}$ and  $L\Sigma_{1}^{-(1)}=L\Sigma_{1}^{-}=I\Pi_{1}^{-}$.

In this Appendix we show that $L\Sigma_{1}^{-(k)}$ is equivalent to the $k$-nested applications of $(LimR)$.
To be precise, let $(LimR)^{(k)}\vdash$ denote the derivability in the $k$-nested applications of $(LimR)$:
$(LimR)^{(0)}\vdash$ is nothing  but  ${\sf EA}\vdash$, and
if $(LimR)^{(k)}\vdash\exists m\forall n\geq m\, h(n+1)\leq h(n)$,
then  $(LimR)^{(k+1)}\vdash\exists m\forall n\geq m \, h(n)=h(m)$.

\begin{theorem}\label{th:klim}
$(LimR)^{(k)}\vdash \varphi \Leftrightarrow L\Sigma_{1}^{-(k)}\vdash\varphi$ for any $\varphi$.
\end{theorem}

This is shown by induction on $k$.
The proof is obtained by a slight modification of proofs in \cite{Lev}.

First consider
\[
(LimR)^{(k)}\vdash L\Sigma_{1}^{-(k)}
.\]
Let $<^{(k)}\, (k\geq 1)$ denote the lexicographic order on $k$-tuples of natural numbers.
Also $\langle x_{1},\ldots,x_{k}\rangle^{(k)}$ denotes a(n elementary recursive) bijective coding of $k$-tuples with its  inverses $(n)_{i}^{(k)}\, (1\leq i\leq k)$. 
In what follows the super scripts $(k)$ are omitted.

Then $L\Sigma_{1}^{-(k)}$ says that if there exists an $x$ satisfying $\varphi(x)\equiv\theta((x)_{1},\ldots,(x)_{k})$,
then there exists a minimal such $x$ with respect to $<^{(k)}$.

We can assume that {\sf EA} proves
\begin{equation}\label{eq:lex}
\exists i[\forall j\neq i(x_{j}=y_{j}) \land x_{i}<y_{i}] \to 
\langle x_{1},\ldots,x_{k}\rangle< \langle y_{1},\ldots,y_{k}\rangle
\end{equation}

Now given a $\Delta^{0}_{0}$-formula $\varphi(x_{1},\ldots,x_{k},x_{k+1})$ without parameters,
we want to show $L\Sigma_{1}^{-(k)}$ with $\theta\equiv \exists x_{k+1}\varphi$.

As  in \cite{Lev} some elementary functions $g_{1}, g, h, h'$ are defined successively as follows.

\[
g_{1}(n)=
\left\{
\begin{array}{ll}
n & \mbox{if } \forall y\leq n \lnot\varphi((y)_{1},\ldots,(y)_{k},(y)_{k+1})  \\
\langle (y)_{1},\ldots,(y)_{k}\rangle & \mbox{otherwise with }
y=\mu y\leq n \varphi((y)_{1},\ldots,(y)_{k},(y)_{k+1})
\end{array}
\right.
\]

\[
g(n)=
\left\{
\begin{array}{ll}
\langle (n)_{1},\ldots,(n)_{k}\rangle & \mbox{if }  \exists u\leq (n)_{k+1} \varphi((n)_{1},\ldots,(n)_{k}, u) \\
g_{1}(n) & \mbox{otherwise}
\end{array}
\right.
\]

$h(0)=g(0)$ and
\[
h(n+1)=
\left\{
\begin{array}{ll}
g(n+1)  & \mbox{if } \forall k,m\leq n(k\neq m \to g(k)\neq g(m)) \\
g(n+1) & \mbox{if } \exists m\leq n (g(n+1)=g(m)) \mbox{ and }
g(n+1)<^{(k)}h(n) \\
h(n) &  \mbox{otherwise}
\end{array}
\right.
\]

Observe that $h(n)\leq\max\{g(m):m\leq n\}$, and hence $h$ is elementary.

\[
h'(x)=
\left\{
\begin{array}{ll}
h(x) & \mbox{ if } \exists n\leq x \varphi((n)_{1},\ldots,(n)_{k}, (n)_{k+1}) \\
0 & \mbox{otherwise}
\end{array}
\right.
\]

Then {\sf EA} proves that $h'$ is eventually decreasing with respect to $<^{(k)}$:
$\exists m\forall n\geq m (h'(n+1)\leq^{(k)}h'(n))$.
Therefore $h'_{1}(n)=(h'(n))_{1}$ is eventually decreasing.
Hence $\exists y_{1}[y_{1}=\lim_{x\to\infty}h_{1}'(x)]$ in $(LimR)^{(1)}$.

This in turn implies that
$\langle (h'(n))_{2},\ldots,(h'(n))_{k}\rangle$ is eventually decreasing with respect to $<^{(k-1)}$.
Therefore $h'_{2}(n)=(h'(n))_{2}$ is eventually decreasing demonstrably in $(LimR)^{(1)}$.
Hence $\exists y_{2}[y_{2}=\lim_{x\to\infty}h_{2}'(x)]$ in $(LimR)^{(2)}$, and so on.
Therefore $\exists y[y=\lim_{x\to\infty}h'(x)]$ in $(LimR)^{(k)}$.

Now assuming $\exists x_{1}\cdots\exists x_{k}\exists x_{k+1} \theta(x_{1},\ldots,x_{k}, x_{k+1})$, 
we see as in \cite{Lev} that 
$y=\lim_{x\to\infty}h'(x)=\lim_{x\to\infty}h(x)$, and the limit $y$ is the minimum of 
$\{\langle x_{1},\ldots,x_{k}\rangle: \exists x_{k+1} \theta(x_{1},\ldots,x_{k}, x_{k+1})\}$
with respect to the lexicographic order $<^{(k)}$ as desired.
\\ \smallskip \\

Next assume by IH that
\[
L\Sigma_{1}^{-(k)}\vdash \exists m\forall n\geq m\, h(n+1)\leq h(n)
.\]
We need to show
\[
L\Sigma_{1}^{-(k+1)}\vdash
\exists m\forall n\geq m \, h(n)=h(m)
.\]
For simplicity consider the case $k=1$, and assume that {\sf EA} proves that
\begin{eqnarray*}
&& \{\exists x_{1}\varphi_{1}(x_{1})\to\exists x_{1}[\varphi_{1}(x_{1}) \land \forall z<x_{1}\lnot\varphi_{1}(z)]\} \land 
\\
&&
\{\exists x_{2}\varphi_{2}(x_{2})\to\exists x_{2}[\varphi_{2}(x_{2}) \land \forall z<x_{2}\lnot\varphi_{2}(z)]\} \\
&&  \to 
\exists m\forall n\geq m\, h(n+1)\leq h(n)
\end{eqnarray*}
Let $\varphi_{i}(x_{i})\equiv \exists y \theta_{i}(x_{i},y)$.

By the Herbrand's Theorem there exists a sequence  of  terms
$m_{0}(a_{1},a_{2},b_{1},b_{2})$, $m_{1}(x_{0},a_{1},a_{2},b_{1},b_{2})$, 
$m_{2}(x_{0},x_{1},a_{1},a_{2},b_{1},b_{2}),\ldots,m_{k}(x_{0},\ldots,x_{k-1},a_{1},a_{2},b_{1},b_{2})$ such that the following disjunction is provable in {\sf EA}:
\begin{eqnarray}
&& \{\exists x_{1}\varphi_{1}(x_{1})\land [\lnot\theta_{1}(a_{1},b_{1}) \lor \exists z<a_{1} \varphi_{1}(z)]\} \lor 
\label{eq:dis} \\
&& \{\exists x_{2}\varphi_{2}(x_{2})\land [\lnot\theta_{2}(a_{2},b_{2}) \lor \exists z<a_{2} \varphi_{2}(z)]\} \lor 
\nonumber \\
&& (x_{0}\geq m_{0}(a_{1},a_{2},b_{1},b_{2}) \to h(x_{0})\geq h(x_{0}+1))\lor 
\nonumber \\
&& (x_{1}\geq m_{1}(x_{0}, a_{1},a_{2},b_{1},b_{2}) \to h(x_{1})\geq h(x_{1}+1))\lor 
\nonumber \\
&& \cdots
\nonumber  \\
&& (x_{k}\geq m_{k}(x_{0}, \ldots,x_{k-1},a_{1},a_{2},b_{1},b_{2}) \to h(x_{k})\geq h(x_{k}+1))
\nonumber 
\end{eqnarray}

First assume $\exists x_{1}\varphi_{1}(x_{1})\land\exists x_{2}\varphi_{2}(x_{2})$, and pick a minimal
$a=\langle a_{1},a_{2}\rangle$ such that $\varphi_{1}(a_{1})\land\varphi_{2}(a_{2})$
and $\forall b<a\lnot[\varphi_{1}((b)_{1})\land\varphi_{2}((b)_{2})]$.
Then from (\ref{eq:lex}) we see that 
$\forall z<a_{i}\lnot\varphi_{i}(z)$ for $i=1,2$.

Now let 
\[
y_{0}=\mu y_{0}[\exists b_{1}, b_{2}\exists x_{0}\geq m_{0}(a_{1},a_{2},b_{1},b_{2})
\{\theta_{1}(a_{1},b_{1}) \land \theta_{2}(a_{2},b_{2}) \land h(x_{0})=y_{0}\}]
\]
by $L\Sigma_{1}^{-(2)}$.

If $\forall x_{0}\geq m_{0}(a_{1},a_{2},b_{1},b_{2})(h(x_{0})\geq h(x_{0}+1))$ for some $b_{1},b_{2}$
with $\theta_{1}(a_{1},b_{1}) \land \theta_{2}(a_{2},b_{2})$, then 
$y_{0}=\lim_{x\to\infty}h(x)$.

Otherwise let 
\begin{eqnarray*}
&& y_{1}=\mu y_{1}[\exists b_{1}, b_{2}\exists x_{0}\geq m_{0}(a_{1},a_{2},b_{1},b_{2})
\exists x_{1}\geq m_{1}(x_{0}, a_{1},a_{2},b_{1},b_{2})
\\
&&
 (\theta_{1}(a_{1},b_{1}) \land \theta_{2}(a_{2},b_{2})  \land h(x_{0})< h(x_{0}+1)
\land h(x_{1})=y_{1})]
\end{eqnarray*}
If $\forall x_{1}\geq m_{1}(x_{0}, a_{1},a_{2},b_{1},b_{2}) (h(x_{1})\geq h(x_{1}+1))$
 for some $b_{1},b_{2}$ with $\theta_{1}(a_{1},b_{1}) \land \theta_{2}(a_{2},b_{2})$, 
 and an $x_{0}\geq m_{0}(a_{1},a_{2},b_{1},b_{2})$, then 
$y_{1}=\lim_{x\to\infty}h(x)$, and so on.

If $\lnot\exists x_{1} \varphi_{1}(x_{1})$, then substitute $0$ for  $a_{1}, b_{1}$ in (\ref{eq:dis}).

\end{document}